\documentclass{article}

\usepackage{enumerate}
\usepackage{amsmath,xspace,amssymb,mathrsfs}

\input xy
\xyoption{all}
\CompileMatrices

\title{Bicat is not triequivalent to Gray}
\author{Stephen Lack%
\thanks{The support of the Australian Research Council and
DETYA is gratefully acknowledged.}
\\School of Computing and Mathematics\\
University of Western Sydney\\
Locked Bag 1797 Penrith South DC NSW 1797\\
Australia\\
email: {\tt s.lack@uws.edu.au}}
\date{}

\newcommand{\A}{{\ensuremath{\mathscr A}}\xspace}
\newcommand{\B}{{\ensuremath{\mathscr B}}\xspace}
\newcommand{\C}{{\ensuremath{\mathscr C}}\xspace}
\newcommand{\Gray}{\textnormal{\bf Gray}\xspace}
\newcommand{\Bicat}{\textnormal{\bf Bicat}\xspace}
\newcommand{\twocat}{\textnormal{\bf 2-Cat}\xspace}
\newcommand{\Z}{{\ensuremath{\mathbb Z}}\xspace}

\newcommand{\st}{\textnormal{\bf st}\xspace}

\newtheorem{theorem}{Theorem}
\newtheorem{lemma}[theorem]{Lemma}
\newtheorem{preremark}[theorem]{Remark}
\newenvironment{remark}{\begin{preremark}\rm}{\end{preremark}}

\newcommand{\proof}{\noindent{\sc Proof:}\xspace}
\def\endproof{~\hfill$\Box$\vskip 10pt}

\begin{document}

\label{firstpage}
\maketitle

\begin{abstract}
\Bicat is the tricategory of bicategories, homomorphisms,
pseudonatural transformations, and modifications. \Gray is the
subtricategory of 2-categories, 2-functors, pseudonatural
transformations, and modifications. We show that these two
tricategories are not triequivalent.
\end{abstract}


Weakening the notion of 2-category by replacing all 
equations between 1-cells by suitably coherent isomorphisms
gives the notion of {\em bicategory} \cite{bicategories}. The analogous 
weakening of a 2-functor is called a {\em homomorphism} of bicategories, 
and the weakening of a 2-natural transformation is a {\em pseudonatural 
transformation}. There are also {\em modifications} between 2-natural
or pseudonatural 
transformations, but this notion does not need to be weakened. The
bicategories, homomorphisms, pseudonatural transformations, and modifications 
form a tricategory (a weak 3-category) called \Bicat. 

The subtricategory of \Bicat containing only the 2-categories as 
objects, and only the 2-functors as 1-cells, but with all 2-cells and
3-cells between them, is called \Gray. As well as being a particular
tricategory, there is another important point of view on \Gray. The
category \twocat of 2-categories and 2-functors is cartesian closed,
but it also has a different symmetric monoidal closed structure \cite{Gray}, 
for which the internal hom $[\A,\B]$ is the 2-category of 2-functors,
pseudonatural transformations, and modifications
between \A and \B. A category enriched over \twocat with respect 
to this closed structure is called a {\em Gray-category}. A Gray-category has
2-categories as hom-objects, so is a 3-dimensional categorical
structure, and it can be seen as a particular sort of tricategory.
The closed structure of \twocat gives it a canonical enrichment over 
itself and the resulting Gray-category is just \Gray. \Gray is also
sometimes used as a name for \twocat with this monoidal structure.

A homomorphism of bicategories $T:\A\to\C$ is called a {\em biequivalence}
if it induces equivalences $T_{A,B}:\A(A,B)\to\C(TA,TB)$ of hom-categories
for all objects $A,B\in\C$ ($T$ is {\em locally an equivalence}), and every 
object $C\in\C$ is equivalent in \C to one of the form $TA$ ($T$ is 
{\em biessentially surjective on objects}). We then write $\A\sim\B$.
Every bicategory is equivalent to a 2-category \cite{MacLane-Pare}. 

A trihomomorphism 
of tricategories $T:\A\to\C$ is called a {\em triequivalence} if it induces
biequivalences $T_{A,B}:\A(A,B)\to\C(TA,TB)$ of hom-bicategories for
all objects $A,B\in\A$ ($T$ is {\em locally a biequivalence}), and every object
$C\in\C$ is biequivalent in \C to one of the form $TA$ ($T$ is 
{\em triessentially surjective on objects}). It is not the case that every 
tricategory is triequivalent to a 3-category, but every tricategory is 
triequivalent to a Gray-category~\cite{GPS}.

Perhaps since a \Gray-category is a category enriched in the monoidal category
\Gray, and a tricategory can be seen as some sort of ``weak
\Bicat-category'', it has been suggested that \Bicat might be
triequivalent to \Gray, and indeed Section~5.6 of \cite{GPS} states
that this is the case. We prove that it is not.
First we prove:

\begin{lemma}
The inclusion $\Gray\to\Bicat$ is not a triequivalence.
\end{lemma}

\proof
If it were then each inclusion $\Gray(\A,\B)\to\Bicat(\A,\B)$ would
be a biequivalence, and so each homomorphism (pseudofunctor) between
2-categories would be pseudonaturally equivalent to a 2-functor. This
is not the case. For example (see \cite[Example~3.1]{qm2cat}), let \A be the 2-category 
with a single object $*$, a single non-identity morphism $f:*\to *$ 
satisfying $f^2=1$, and no non-identity 2-cells (the group of order 2
seen as a one-object
2-category); and let $\B$ be the 2-category with a single object
$*$, a morphism $n:*\to *$ for each integer $n$, composed via addition,
and an isomorphism $n\cong m$ if and only if $n-m$ is even (the 
``pseudo-quotient of \Z by $2\Z$''). There is
a homomorphism $\A\to\B$ sending $f$ to $1$; but the only 2-functor
$\A\to\B$ sends $f$ to 0, so this homomorphism is not pseudonaturally
equivalent to a 2-functor.
\endproof

\begin{theorem}
\Gray is not triequivalent to \Bicat.
\end{theorem}

\proof
Suppose there were a triequivalence $\Phi:\Gray\to\Bicat$.
We show that $\Phi$ would 
be biequivalent to the inclusion, so that the inclusion itself would
be a triequivalence; but by the lemma this is impossible.

The terminal 2-category 1 is a terminal object in \Gray, so must be sent to a 
``triterminal object'' $\Phi 1$ in \Bicat; in other words,
$\Bicat(\B,\Phi 1)$ must be biequivalent to 1 for any bicategory 
\B. For any 2-category \A, we have biequivalences
$$\A\sim\Gray(1,\A)\sim\Bicat(\Phi 1,\Phi\A)\sim\Bicat(1,\Phi\A)\sim\Phi\A$$
where the first is the isomorphism coming from the monoidal structure on \Gray,
the second is the biequivalence on hom-bicategories given by $\Phi$, the
third is given by composition with the biequivalence $\Phi 1\sim 1$, and the 
last is a special case of the biequivalence $\Bicat(1,\B)\sim\B$ for any 
bicategory, given by evaluation at the unique object $*$ of $1$. All of these 
biequivalences are ``natural'' in a suitably weak tricategorical
sense, and so $\Phi$ is indeed biequivalent to the inclusion. 
\endproof

\begin{remark}
The most suitable weak tricategorical transformation is called 
a tritransformation. The axioms are rather daunting, but really
the coherence conditions are not needed here. We only need the obvious
fact that for any 2-functor $T:\A\to\B$, the square
$$\xymatrix{
\A \ar[d]_{T} \ar[r]^{\sim} & \Phi\A \ar[d]^{\Phi T} \\
\B \ar[r]_{\sim} & \Phi\B }$$
commutes up to equivalence.
\end{remark}

The fact that every bicategory is biequivalent to a 2-category is
precisely the statement that the 
inclusion $\Gray\to\Bicat$ is triessentially surjective on objects,
but as we saw in the lemma, it is not locally a biequivalence.
On the other hand Gordon, Power, and Street construct in \cite{GPS} a 
trihomomorphism $\st:\Bicat\to\Gray$ which is locally a biequivalence
(it induces a biequivalence on the hom-bicategories). They do this 
by appeal to their Section 3.6, but this does not imply that $\st$ 
is a triequivalence, as they claim, and by our theorem it cannot be one.
In fact Section~5.6 is not used in the proof of the main theorem of
\cite{GPS}, it is only used to construct the tricategory \Bicat
itself, and this does not need $\st$ to be a triequivalence.

By the coherence result of \cite{GPS}, \Bicat is triequivalent
to {\em some} Gray-category; and by the fact that \st\ is locally a
biequivalence, \Bicat is triequivalent to a full sub-Gray-category of \Gray, 
but it is not triequivalent to \Gray itself.


\begin{thebibliography}{9}

\bibitem{bicategories}
Jean B{\'e}nabou.
\newblock Introduction to bicategories.
\newblock In {\em Reports of the Midwest Category Seminar}, pages 1--77.
  Springer, Berlin, 1967.

\bibitem{GPS}
R.~Gordon, A.~J. Power, and Ross Street.
\newblock Coherence for tricategories.
\newblock {\em Mem. Amer. Math. Soc.}, 117(558):vi+81, 1995.

\bibitem{Gray}
John~W. Gray.
\newblock {\em Formal category theory: adjointness for {$2$}-categories}.
\newblock Springer-Verlag, Berlin, 1974.

\bibitem{qm2cat}
Stephen Lack.
\newblock A {Q}uillen model structure for 2-categories.
\newblock {\em $K$-Theory}, 26(2):171--205, 2002.

\bibitem{MacLane-Pare}
Saunders Mac~Lane and Robert Par{\'e}.
\newblock Coherence for bicategories and indexed categories.
\newblock {\em J. Pure Appl. Algebra}, 37(1):59--80, 1985.
\end{thebibliography}

\end{document}